\documentclass[11pt,a4paper]{article}

\usepackage{amssymb}
\usepackage{graphics}
\usepackage{amsmath}
\usepackage{a4wide}
\usepackage[ruled]{algorithm2e}

\newcommand{\lbn}{\ensuremath{L_B(n)}}
\newcommand{\lbnp}{\ensuremath{L_B(n+1)}}
\newcommand{\lbi}{\ensuremath{L_B(\infty)}}
\newcommand{\tlbi}{\ensuremath{\widetilde{L_B(\infty)}}}
\newcommand{\tb}{\ensuremath{T_B(\infty)}}

\newcommand{\ie}{\emph{i.e.}}

\newcommand{\tr}[1]{\ensuremath{\stackrel{#1}{\longrightarrow}}}
\newcommand{\fb}[1]{\ensuremath{{\downarrow_{#1}}}}
\newcommand{\inc}[2]{\ensuremath{#1^\fb{#2}}}

\newcommand{\lgr}[2]{\raisebox{-2.5ex}{%
$\stackrel{\underleftrightarrow{\mbox{\normalsize $#1$}}}{\mbox{\tiny $^{#2}$}}$}}

\newcommand{\includeferrers}[1]{\scalebox{0.2}{\includegraphics{#1.eps}}}

\newcommand{\fig}[2]{
\begin{figure}[!h]
\centering
\includegraphics{#1.eps}
\caption{#2}
\label{fig_#1}
\end{figure}
}

\begin{document}

\newenvironment{proof}{\noindent{\bf Proof:}}{\ \hfill\ $\square$}

\newtheorem{theorem}{Theorem}
\newtheorem{corol}{Corollary}
\newtheorem{lemme}{Lemma}
\newtheorem{definition}{Definition}
\newtheorem{prop}{Proposition}

\begin{center}

{\bf \LARGE The lattice of integer partitions\\\smallskip and its infinite extension}

\bigskip

{ Matthieu Latapy\ \ \ \ Thi Ha Duong Phan}

\bigskip
\bigskip

\begin{minipage}{.9\textwidth}

\centerline{\bf Abstract.}

In this paper, we use a simple discrete
dynamical model to study integer partitions and their lattice. The
set of reachable configurations of the model, with the order induced
by the transition rule defined on it, is the lattice of all partitions
of an integer, equipped with a dominance ordering. We first
explain how this lattice can be constructed by an algorithm in linear
time with respect to its size by showing that it has a 
self-similar structure. Then, we define a natural extension of the
model to infinity, which we compare with the Young lattice. Using a
self-similar tree, we obtain an encoding of the obtained
lattice which makes it possible to enumerate easily and efficiently
all the partitions of a given integer. This approach also gives a
recursive formula for the number of partitions of an integer, and some
informations on special sets of partitions, such as length bounded
partitions.

\medskip


\end{minipage}

\end{center}

\bigskip

\section{Preliminaries}

A \emph{partially ordered set} (or \emph{poset}) is a set $P$ with a
reflexive ($x \le x$), transitive ($x \le y$ and $y \le z$ implies
$x \le z$) and antisymmetric ($x \le y$ and $y \le x$ implies $x=y$)
binary relation $\le$. A {\em lattice} is a partially ordered set
such that any two elements $a$ and $b$ have a least upper bound,
called {\em supremum} of $a$ and $b$ and denoted by $sup(a,b)$, and
a greatest lower bound, called {\em infimum} of $a$ and $b$ and
denoted by $inf(a,b)$. The element $sup(a,b)$ is the smallest
element among the elements greater than both $a$ and $b$. The
element $inf(a,b)$ is defined dually. A subset $L'$ of a lattice $L$ is called a {\em sublattice} of $L$ if for any two elements $a$ and $b$ of $L'$, the $sup(a,b)$ and $inf(a,b)$ are also elements of $L'$. Lattices are strongly
structured sets, and many general results, for instance efficient
encodings and algorithms, are known about them. For more details,
see for instance \cite{DP90}.

A \emph{partition} is an integer sequence $a = (a_1,a_2,\dots,a_k)$
such that $a_1 \geq a_2 \geq \ldots \geq a_k >0$ (by convention, $a_j =0$ for all $j>k$). We say that $a$ is a partition of $n$ if $\sum_{i=1}^{i=k} a_i =n$.
The Ferrers diagram of a partition $a = (a_1,a_2,\dots,a_k)$ is a
drawing of $a$ on $k$ adjacent columns such that the $i$-th column
is a pile of $a_i$ stacked squares, which we will call \emph{grains}
because of the sand piles dynamics we will consider over them.
For instance, $p=(4,3,3,2)$
and $q=(6,2,1,1,1,1)$ are two partitions of $n=12$,
and their Ferrers diagrams are
\includeferrers{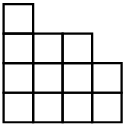}
and
\includeferrers{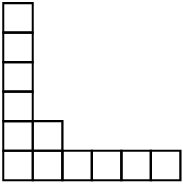}
respectively.

The \emph{dominance ordering} is defined in the following way \cite{Bry73}.
Consider two partitions of the integer $n$: $a = (a_1,a_2,\dots,a_k)$ and
$b = (b_1,b_2,\dots,b_l)$. Then
$$a \ge b \mbox{ if and only if }
\sum_{i=1}^j a_i \ge \sum_{i=1}^j b_i \mbox{ for all $j \geq 1$.}
$$

From \cite{Bry73}, it is known that the set of all partitions of an
integer $n$ with the dominance ordering is a lattice, denoted
by \lbn. In his paper, Brylawski proposed a dynamical approach to study
this lattice. We will introduce some notations to explain
it intuitively. For more details about integer partitions, we refer
to \cite{And76}.

Let $a = (a_1,\dots a_k)$ be a partition. The
\emph{height difference of $a$ at $i$},
denoted by $d_i(a)$, is the integer $a_i - a_{i+1}$.
We say that the partition $a$ has a \emph{cliff} at $i$ if $d_i(a) \ge 2$.
We say that $a$ has a \emph{slippery plateau} at $i$ if there exists $\ell>i$
such that $d_j(a) = 0$ for all
$i \le j < \ell$ and $d_{\ell}(a) = 1$. The integer $\ell-i$ is then called the
\emph{length} of the slippery plateau at $i$.
Likewise, $a$ has a \emph{non-slippery plateau} at $i$ if $d_j(a) = 0$
for all $i \le j < \ell$ and it has a cliff at $\ell$. The integer $\ell-i$ is
called the \emph{length} of the non-slippery plateau at $i$.
The partition $a$ has a \emph{slippery step} at $i$ if the sequence
defined by $a' = (a_1, \dots, a_i-1, \dots, a_k)$ is a partition with a
slippery plateau at $i$. Likewise, $a$ has a \emph{non-slippery step}
at $i$ if $a'$ is a partition with a non-slippery plateau at $i$.
See Figure~\ref{fig_defs} for some illustrations.

\fig{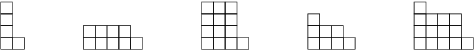}{From left to right: a cliff, a slippery plateau of length $3$,
a non-slippery plateau of length $2$, a slippery step of length $2$ and
a non-slippery step of length $3$.}

\fig{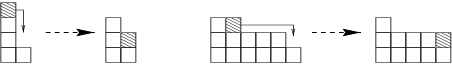}{The two evolution rules of the dynamical model}

Consider now the partition $a=(a_1,a_2,\dots,a_k)$.
Brylawski defined the following two evolution rules:
one grain can fall from column $i$ to column $i+1$ if $a$ has a cliff at $i$,
and one grain can slip from column $i$ to column $i+l+1$ if $a$ has a
slippery step of length $l$ at $i$. See Figure~\ref{fig_regles}.

Such a fall
or a slip is called a \emph{transition} of the model and is denoted by
$a \tr{i} b$ where $i$ is the column from which the grain falls or slips; we also denote this by $b = a \tr{i}$.
If one starts from the partition $(n)$ and
iterates this operation,
one obtains all the partitions of $n$, and the dominance ordering is
nothing but the reflexive and transitive closure of the relation
induced by the transition rule \cite{Bry73}.
See Figure~\ref{fig_diag} for illustrations with $n=7$ and $n=8$.

\fig{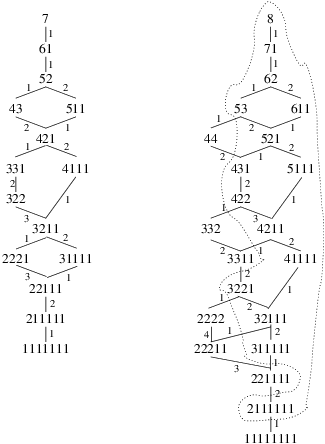}{Diagrams of the lattices $L_B(n)$ for $n=7$ and $n=8$. As we will
see, the set $L_B(7)$ is isomorphic to a sublattice of $L_B(8)$. On the
diagram of $L_B(8)$, we included in a dotted line this sublattice.}

Let us recall that one can consider Brylawski's model as a generalization
of the so-called Sand Pile Model (SPM), which consists of the first
evolution rule only. The SPM was studied in many areas: from physics
point of view \cite{BTW87}, combinatorics considerations
\cite{ALSSTW89,GK93}, and dynamical model theory
\cite{GMP02a,GMP02b,LP01}. Moreover an infinite extension of this
model was studied in \cite{LMMP01}.

We will now study the structure of the lattice of the partitions
of an integer $n$ and we will show its self-similarity by giving
a method to construct \lbnp\ from \lbn. Then,
we will define an infinite extension of these lattices: the
lattice \lbi\ of all the partitions of any integer (\ie\ all finite
non-increasing sequences of positive integers). 
We will compare this lattice with the Young lattice, which also contains
all the partitions of any integer, but ordered in a different way.
Finally, we will construct an infinite tree based on the
construction process described at the beginning of the paper.
This tree will make it possible to give a simple and efficient algorithm
to enumerate all the partitions of a given integer. It also has
a self-similar structure, from which we will obtain a 
recursive formula
for the number of partitions of an integer $n$ and some results
about certain classes of partitions.

Before entering the core of the topic, we need one more notation.
If the $k$-tuple \mbox{$a=(a_1,a_2,\dots,a_k)$} is a partition,
then the $k$-tuple
$(a_1,a_2,\dots,a_{i-1},a_i+1,a_{i+1},\dots,a_k)$ is denoted by \inc{a}{i}.
In other words, \inc{a}{i} is obtained from $a$ by adding one grain
on its $i$-th column.
Notice that the $k$-tuple obtained this way is not
necessarily a partition. If $S$ is a set of partitions,
then \inc{S}{i}\ denotes the set $\lbrace \inc{a}{i} | a \in S \rbrace$.
Finally, we denote by
$dirreach(a)$
the set of configurations directly reachable from
$a$, \ie\
the set $\lbrace b\ |\ a \tr{i} b \mbox{ for some } i \rbrace$.
Notice that in the context of dynamical model theory, those
elements are called \emph{the immediate successors} of $a$. However, since
we are concerned here with order theory, we cannot use this term,
which takes another meaning in this context.

\section{From \lbn\ to \lbnp}

\label{sec_deb}

In this section, our aim is to construct \lbnp\ from \lbn, viewed as
the graph induced by the dynamical model, with the edges labeled by
the number of the column from which the grain falls or slips, as
shown in Figure~\ref{fig_diag}. We will call \emph{construction of a
lattice} the computation of this labeled graph. We first show that
\inc{\lbn}{1}\ is a sublattice of \lbnp. For instance, in
Figure~\ref{fig_diag} we included in a dotted line \inc{L_B(7)}{1}\
within $L_B(8)$. This remark allows us to start the construction of
\lbnp\ from \lbn\ by computing \inc{\lbn}{1}\ and then adding the
missing elements of \lbnp. After characterizing those elements that
must be added, we obtain a simple and efficient method to achieve
the construction of \lbnp\ from \lbn.

\begin{prop}
\label{prop_sub}
\inc{\lbn}{1}\ is a sublattice of \lbnp.
\end{prop}
\begin{proof}
We must show that for any two elements $a$ and $b$ of \lbn, $inf(\inc{a}{1},\inc{b}{1})$ and $sup(\inc{a}{1},\inc{b}{1})$ are in $\inc{\lbn}{1}$.
 Let us first consider the element $c=inf(a,b)$. It is clear that $\inc{c}{1}$ is in $\inc{\lbn}{1}$, and we will show that $\inc{c}{1}$ is equal to 
$inf(\inc{a}{1},\inc{b}{1})$. This statement comes directly from Brylawski's result on dominance ordering \cite{Bry73}:
$$
inf(a,b) = c\mbox{ if and only if, for all $j \geq 1$, one has }\sum_{i=1}^j c_i =
min(\sum_{i=1}^j a_i,\sum_{i=1}^j b_i).
$$

Let us consider now $e = sup(a,b)$. We will show that $\inc{e}{1}$ is equal to
$d = sup(\inc{a}{1},\inc{b}{1})$.
We have $e \ge a$ and $e \ge b$, therefore
$\inc{e}{1} \ge \inc{a}{1}$ and $\inc{e}{1} \ge
\inc{b}{1}$. This implies that $\inc{e}{1} \ge d$.
To show that $d \ge \inc{e}{1}$,
let us begin by showing that $d_1 = c_1 + 1$. We can suppose that $a_1 \ge b_1$.
The partition $(a_1, a_1, a_1-1, a_1-2, \dots )$ is greater than $a$ and $b$,
and so it is greater than or equal to $e$. Moreover, $e \ge a$ implies $e_1 \ge a_1$
and so $e_1 = a_1$. Since $\inc{a}{1} \le d \le \inc{e}{1}$,
we then have $d_1 = a_1+1 = e_1+1$.
Let $f=(d_1-1, d_2, d_3,\dots)$. Since $d \le \inc{e}{1}$ and $e_1 = a_1$,
$f$ is a partition: $d_1-1 \ge d_2$. Moreover, $d \ge \inc{a}{1}$
and $d \ge \inc{b}{1}$, and so $f \ge a$ and $f \ge b$.
This implies that $f \ge sup(a,b) = c$ and that
$d \ge \inc{e}{1}$, which ends the proof.
\end{proof}

This result shows that one can construct the lattice \lbnp\ from
\lbn\ as follows. The first step of this construction is to
construct the set \inc{\lbn}{1}\ by adding one grain to the first
column of each element of \lbn. Then, one has to add the missing
elements and their transitions. Therefore, we will now consider the
consequences of the addition of one grain on the first column of a
partition, depending on its structure.

It is clear that, for $n\geq 1$, $\lbn = C(n) \bigsqcup S(n) \bigsqcup NS(n) \bigsqcup NP(n) \bigsqcup_{n> \ell \ge 1} P_{\ell}(n)$, where
 $C(n), S(n), NS(n), P_{\ell}(n)$ and $NP(n)$ are respectively the set of partitions of $n$ with a cliff at $1$, with a slippery step at $1$,  with a non-slippery step at $1$, with a slippery plateau at $1$,
and with a non-slippery plateau of at $1$, and where $\sqcup$ denotes the disjoint union.

\begin{prop}
\label{prop_prelim}
Let $a$ be a partition. Then, we have:
\begin{enumerate}
\item if $a \in C(n)$ or $a \in NP(n)$ then\\
 \centerline{$dirreach(\inc{a}{1}) = \inc{dirreach(a)}{1}$;}
\item if $a \in  P_{\ell}(n)$ then 
 $\inc{a}{1} \tr{1} \inc{a}{\ell +2}$ and\\
 \centerline{$dirreach(\inc{a}{1}) = \inc{dirreach(a)}{1}
                       \cup \lbrace \inc{a}{\ell +2} \rbrace$;}
\item if $a \in S(n)$ 
and  $b$ is such that $a \tr{1} b$, then we have
 $\inc{a}{1} \tr{1} \inc{a}{2} \tr{2} \inc{b}{1}$ and\\
 \centerline{$dirreach(\inc{a}{1}) = \inc{(dirreach(a)
                                      \setminus \lbrace b \rbrace)}{1}
                   \cup \lbrace \inc{a}{2} \rbrace$;}
\item if $a \in  NS(n)$ then 
 $\inc{a}{1} \tr{1} \inc{a}{2}$ and\\
 \centerline{$dirreach(\inc{a}{1}) = \inc{dirreach(a)}{1}
                       \cup \lbrace \inc{a}{2} \rbrace$.}                  
\end{enumerate}
\end{prop}
\begin{proof}
It is obvious that the right hand side of each of these equations is a 
subset of its correspond left hand side, thus it is sufficient to prove the converse. 
So, let us consider an element $c = \inc{a}{1} \tr{i}$ in  $dirreach(\inc{a}{1})$.

\begin{enumerate}
\item
If $i \neq 1$ then $c = \inc{(a \tr{i})}{1} \in \inc{dirreach(a)}{1}$. 
If $a \in NP(n)$ then there is no transition at the first column of $\inc{a}{1}$.
If $a \in C(n)$, and if $c= \inc{a}{1} \tr{1}$ then $c$ is also equal to $\inc{(a \tr{1})}{1}$.
\item
If $i \neq 1$ then $c = \inc{(a \tr{i})}{1} \in \inc{dirreach(a)}{1}$. 
Otherwise if $i=1$, it is clear that $c = \inc{a}{\ell +2}$.
\item
If $i \neq 1$ then $c = \inc{(a \tr{i})}{1} \in \inc{dirreach(a)}{1}$. 
However, the element $b=a \tr{1}$ is obtained directly from $a$, but $\inc{b}{1}$ is not obtained directly from $\inc{a}{1}$. So we have:
$c \in \inc{(dirreach(a) \setminus \lbrace b \rbrace)}{1}$.
Otherwise if $i=1$, it is clear that $c = \inc{a}{2}$.
\item 
If $i \neq 1$ then $c = \inc{(a \tr{i})}{1} \in \inc{dirreach(a)}{1}$. 
Otherwise if $i=1$, it is clear that $c = \inc{a}{2}$.
\end{enumerate}
\end{proof}

In the following theorem, we will represent the set \lbnp\ as a disjoint union.
In addition, its proof will give the transitions between the elements of this set,
which will complete the construction of the lattice \lbnp.

\begin{theorem}
\label{th_complet}
For all $n \geq 1$, we have:
$$L_B(n+1) = \inc{L_B(n)}{1} \bigsqcup \inc{S(n)}{2} \bigsqcup \inc{NS(n)}{2} 
\bigsqcup_{n > \ell \ge 1} \inc{P_{\ell}(n)}{\ell +2}.$$
\end{theorem}
\begin{proof}
First of all, it is easy to check that this union is a disjoint union.

Let us recall the strategy of the construction of \lbnp. First, \inc{\lbn}{1} is a sublattice of \lbnp; we then add to \inc{\lbn}{1} all directly reachable elements (and transitions) from this set to obtain a new set $L^1$; finally, we add to $L^1$ all  directly reachable elements (and transitions) from $L^1$ to obtain a new set $L^2$, and so on. The key idea of this theorem is to show that the set $L^1$ is already equal to \lbnp, or, equivalently, that all directly reachable elements from $L^1$ are elements of $L^1$. 

In Proposition 2, we have shown that $L^1$ is represented as the disjoint union in the right hand side of the claim. 
So let $g \in L^1$, and let $h = g \tr{i}$ be a directly reachable element of $g$; we shall prove that $h \in L^1$.
Several cases are possible.
\begin{trivlist}
\item[$\ \bullet$] $g= \inc{a}{1}$ with $a \in \lbn$. From Proposition 2, all directly reachable elements from $g$ are in $L^1$.

\item[$\ \bullet$] $g = \inc{a}{2}$ with $a \in S(n)$. The transition $a \tr{1} c$ is possible in \lbn. All transitions $g = \inc{a}{2} \tr{i} h = \inc{b}{2}$ are the same as transition $a \tr{i} b$, except the transition $a \tr{1} c$. 
Moreover, it is clear that, if $b$ belongs to $S(n)$, then $h$ belongs to $\inc{S(n)}{2}$. Regarding $c$, we have the transition $\inc{a}{2}  \tr{2} \inc{c}{1}$, and $h = \inc{c}{1}$ belongs to $\inc{L_B(n)}{1}$.

\item[$\ \bullet$] $g= \inc{a}{2}$ with $a \in NS(n)$. All transitions $g \tr{i} h$ are the same as transition $a \tr{i} b$, and $h = \inc{b}{2}$. If $i >2$ then $b$ belongs to $NS(n)$, and then $h$ belongs to $\inc{NS(n)}{2}$. Otherwise, $i$ can be $2$ in the case where $a$ has a non-slippery step of length 1 at 1. In this case, $b$ has a cliff at 1, and $\inc{b}{2} = \inc{c}{1}$ with $c \in L_B(n)$. Hence $h = \inc{b}{2} \in \inc{L_B(n)}{1}$.

\item[$\ \bullet$] $g= \inc{a}{\ell+2}$ with $a \in P_{\ell}(n)$. This case requies more attention. We distinguish three subcases:
\begin{enumerate}
\item $a$ has a cliff at $\ell+1$. Then all transitions $g \tr{i} h$ are the same as transition $a \tr{i} b$, and $h = \inc{b}{\ell+2}$. Moreover, $b$ is an element of $P_{\ell}(n)$ so $h \in \inc{P_{\ell}}{\ell+2}$.
\item $a$ has a non-slippery step at $\ell+1$. Then all transitions $g \tr{i} h$ are the same as transition $a \tr{i} b$, and $h = \inc{b}{\ell+2}$. Moreover, $b$ is an element of $P_{\ell}(n)$ so $h \in \inc{P_{\ell}}{\ell+2}$.
\item $a$ has a slippery step at $\ell+1$. The transition $a \tr{\ell+1} c$ is possible. 
All transitions $g \tr{i} h$ are the same as transition $a \tr{i} b$, except the transition $a \tr{\ell+1} c$. It is easy to check that if $b$ is in $P_{\ell}$, then $h = \inc{b}{\ell+2} \in \inc{P_{\ell}}{\ell+2}$.  Regarding $c$,  
we have the transition $\inc{a}{\ell+2} \tr{\ell+2} \inc{c}{\ell+1}$. Moreover $c$ is an element of $P_{\ell-1}$, so  $h = \inc{c}{\ell+1} \in \inc{P_{\ell -1}}{\ell+1} \in L^1$. This completes the proof.
\end{enumerate}                           
\end{trivlist}
\end{proof}

This result makes it possible to write an algorithm which
constructs the lattice \lbnp\ from \lbn\
in linear time with
respect to the number of added elements and transitions.
Notice that we can obtain \lbn\
for an arbitrary integer $n$ by starting from $L_B(0)$ and iterating this
algorithm,
and so we have an algorithm that constructs \lbn\ in linear time with
respect to its size.

\section{The infinite lattice $L_B(\infty)$}

We will now define \lbi\ as the set of all
configurations reachable from $(\infty)$ (this is the configuration where the
first column contains infinitely many grains and all the other columns
contain no grain).
Therefore, each element $a$ of \lbi\ has the form
$(\infty,a_2,a_3,\dots,a_k)$.
As in the previous section, the dominance ordering
on \lbi\ (when the first component
is ignored) is equivalent to the order induced by
the dynamical model.
The first partitions in \lbi\ are given in
Figure~\ref{fig_lb_infi} along with their covering relations (the first
component, equal to $\infty$, is not represented
on this diagram).

\fig{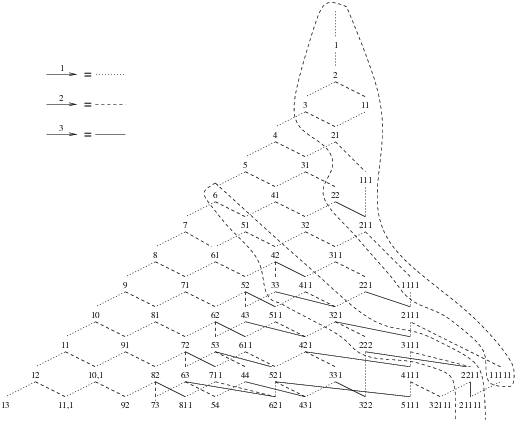}{The first elements and transitions of \lbi. As shown on this
figure for $n=6$, we will discuss two ways to find parts of \lbi\ isomorphic to
\lbn\ for any $n$.}

It is easy to observe that we have a characterization of the order similar to
the one given in \cite{Bry73} for the finite case:
let $a$ and $b$ be two elements of \lbi, $a$ being of length $p$ and
$b$ being of length $q$. Then,
$$
a \ge_{\lbi} b \mbox{ if and only if for all $j$ between $2$ and
$\max(p,q)$, } \sum_{i\ge j} a_i \le \sum_{i\ge j} b_i.
$$

We will start this section by proving that \lbi\ is a lattice and by
giving a formula for the infimum in \lbi. After this, we will show that,
for any $n$,
there are two different ways to find sublattices of \lbi\ isomorphic to
\lbn. We will also give a way to construct some other special sublattices
of \lbi, using its self-similarity. Finally, we will compare \lbi\ with
the Young lattice.

\begin{theorem} \label{th_lbi}
The set \lbi\ is a lattice. Moreover, if
\mbox{$a=(\infty,a_2,\dots,a_k)$} and \mbox{$b=(\infty,b_2,\dots, b_l)$}
are two elements of \lbi, then
$\inf_{\lbi}(a,b)=c$ in \lbi, where $c$ is defined by:
$$
c_i = max(\sum_{j\geq i} a_j, \sum_{j \geq i} b_j) - \sum_{j >i} c_j\hskip
0.3cm \mbox{ for all } i\mbox{ such that } 2 \leq i\leq max(k,l).
$$
\end{theorem}
\begin{proof}
We shall prove that $c$ is an element
of $L_B(\infty)$ and that $c$ is equal to $inf_{\lbi}(a,b)$.
Let $n = 2(\sum_{i\geq 2}a_i + \sum_{i \geq 2}b_i)$.
Let $a'=(n-\sum_{i\geq 2}a_i,a_2,\ldots,a_k)$,
$b'=(n-\sum_{i\geq 2}b_i,b_2,\ldots,b_l)$ and
$c'=(n-\sum_{i\geq 2}c_i,c_2,\ldots,c_{max(k,l)})$.
It is then obvious that $a'$ and $b'$ are two partitions of $n$ and that $c'$
is the infimum
of $a'$ and $b'$ by the dominance ordering in \lbn. Therefore, $c'$ is a
decreasing sequence, and so $c$ is an element of $L_B(\infty)$.
Moreover, according to the definition of $\geq_{\lbi}$, $c$ is
the maximal element of $L_B(\infty)$ which is smaller than $a$ and $b$,
and so $c=inf_{\lbi}(a,b)$.

By definition, $L_B(\infty)$ has a maximal element. Since it is
closed for the
infimum, \lbi\ is a lattice.
\end{proof}

Let us consider now the injective map
$$
\begin{array}{lclc}
\pi : & \lbn                   & \longrightarrow & \lbi\\
     & a=(a_1, a_2,\dots,a_k) & \longrightarrow & \bar{a} =
      (\infty,a_2,\dots,a_k)
\end{array}
$$
One can apply a proof similar to the one of Proposition~\ref{prop_sub} to show that
$$inf_{\lbi}(\pi(a),\pi(b))=\pi(inf_{L_B(n)}(a,b)) \mbox{ and }$$
$$sup_{\lbi}(\pi(a),\pi(b))=\pi(sup_{L_B(n)}(a,b)).$$
This implies that $\pi$ is a lattice embedding.

Let $\overline{\lbn} = \pi(\lbn)$. We know that $\overline{\lbn}$ is a sublattice of \lbi\, and from
Proposition~\ref{prop_sub}, \inc{\lbn}{1} is a sublattice of \lbnp, therefore, since
$\overline{\inc{\lbn}{1}} = \overline{\lbn}$, we have an increasing sequence
of sublattices:
$$
\overline{L_B(0)} \leq \overline{L_B(1)}\leq  \dots \leq \overline{L_B(n)}
\leq \overline{L_B(n+1)} \leq \dots \leq \lbi
$$
where $\le$ denotes the sublattice relation.

We can say more about this increasing sequence of lattices.
Let $a=(\infty,a_2,a_3,\dots,a_k)$ be an element of \lbi.
If one takes $a_1 = a_2 + 1$
and $n = \sum_{i=1}^{k} a_i$, then the partition
$a' = (a_1,a_2,\dots,a_k)$ is an
element of \lbn. Since $a=\pi(a')$, this implies that $a$ is an element of
$\overline{\lbn}$.
Conversly, any element of \lbi\ is of the form $a = (\infty,a_2,\dots,a_k)$.
Therefore, $a'=(a_2,\dots,a_k)$ is a decreasing sequence, and if we
put $n = \sum_{i\ge 2} a_i$ then $a' \in \lbn$, \ie\ $a \in \overline{\lbn}$.
Finally, we have:
$$
\bigcup_{n\ge 0} \overline{\lbn} = \lbi
$$
Therefore, \lbi\ can be viewed as the limit of \lbn\ when $n$ grows to
infinity.
On the other hand, we will show that \lbi\ can be represented as an disjoint union of \lbn\ for all $n$.

Let us define the set 
$$ \tlbi = \bigsqcup_{n\ge0}\lbn. $$
We define the following
relations over \tlbi. Let $a \in L_B(m)$ and $b \in L_B(n)$. We have
$a \tr{i} b$ in \tlbi\ if and only if one of the following applies:
$n=m$ and $a \tr{i} b$ in \lbn, or $i=0$, $n = m+1$ and
$b=\inc{a}{1}$. In other terms, the elements of \lbn\ are linked to
each other as usual, and each element $a$ of \lbn\ is linked to
$\inc{a}{1} \in \lbnp$ by an edge labeled by $0$.

From this, one can introduce an order on the set \tlbi\ in the usual
sense, by defining it as the reflexive and transitive closure of
this relation. We now show that \lbi\ is isomorphic to \tlbi,
and so that \tlbi\ is a lattice.
\begin{lemme}
\label{th_chi}
The map $\chi$ defined by:
$$
\chi : \tlbi \longrightarrow \lbi
$$
$$
a=(a_1,a_2,\ldots,a_k) \mapsto \chi(a) = (\infty,a_1,a_2,\ldots,a_k)
$$
is a lattice isomorphism.\\
 Moreover,
$a \tr{i} b$ in \tlbi\ if and only if $\chi(a) \tr{i+1}
\chi(b)$ in \lbi.
\end{lemme}
\begin{proof}
$\chi$ is clearly bijective. Moreover, it is clear from the
definitions that for all $a$ and $b$ in
\tlbi, $a \tr{i} b$ if and only if $\chi(a) \tr{i+1} \chi(b)$.
Therefore, $\chi$ is an order isomorphism. Since \lbi\ is a lattice,
this implies that $\chi$ is a lattice isomorphism.
\end{proof}

This lemma means that \tlbi\ is nothing but \lbi\ when one removes
the first component (always equal to $\infty$) of each element of \lbi\
and decreases the label of each edge by $1$. We will now see that \lbn\
is a sublattice of \tlbi\ for all $n$, which gives another way to
find a part of \lbi\ isomorphic to \lbn.

\begin{theorem}
For all integer $n \geq 1$, \lbn\ is a sublattice of \tlbi.
\label{th_retrouve2}
\end{theorem}
\begin{proof}
Let $a$ and $b$ be two elements of \lbn, we shall prove that $inf_{\tlbi}(a,b)$
and $sup_{\tlbi}(a,b)$ belong to \lbn. Let $c$ be $inf_{\lbn}(a,b)$ and
$c'$ be $inf_{\tlbi}(a,b)$. We have, $a \geq c' \geq_{\tlbi} c$,
which means that
$\sum_{i \geq 1}a_i \leq \sum_{i \geq 1}c'_i \leq \sum_{i \geq 1}c_i$,
and so $\sum_{i \geq 1}c'_i=n$. This implies that $c'$ belongs to \lbn,
and we obtain $c'=c$. The proof for the supremum is similar.
\end{proof}

To finish this section, we will discuss the relations between the
infinite lattice \tlbi\ and the famous Young lattice. These two
infinite lattices contain exactly the same elements (all the
partitions of all the integers), but ordered in a different way: $a
\le b$ in the Young lattice if for all $i$ we have $a_i \le b_i$. In
other words, the order over the partitions is the componentwise
order. This order induces a (distributive) lattice structure over
the set of all the integer partitions. It has been widely studied;
see for instance \cite{Sta99,Gue93}. It can also be viewed as the set
of partitions obtained from the empty one, $()$, and by iterating
the following evolution rule: $a \tr{i} b$ if $b$ is a partition
obtained from the partition $a$ by increasing its $i$-th component.
This implies directly that the lattice can be decomposed into levels
(the $i$-th level contains the partitions obtained after $i$
applications of the evolution rule), and that level $i$ contains
exactly the partitions of $n$, \ie\ the elements of \lbn. Notice
moreover that these elements are not comparable in the Young lattice
therefore the order in \tlbi\ and the one in the Young lattice are
very different. However, they are in relation according to the
following result:

\begin{prop}\cite{Lat02}
The map $\pi$ from \tlbi\ into the Young lattice
such that $\pi(a)_i$ is equal to $\sum_{j\ge i}a_j$
is an order embedding which preserves the infimum.
\end{prop}
\begin{proof}
Let $a$ and $b$ be two elements of \tlbi. We must show that $\pi(a)$
and $\pi(b)$ belong to the Young lattice, that $a \ge b$ in \tlbi\ is
equivalent to $\pi(a) \ge \pi(b)$ in the Young lattice and that
$inf(\pi(a),\pi(b))$ in the Young lattice is equal to $\pi(inf(a,b))$
in \tlbi.
The first two points are easy:
$\pi(x)$ is obviously a decreasing sequence of integers for any $x$, and
the order is preserved.
Now, let $c = inf(a,b)$. Then,
$$
\begin{array}{lclr}
\pi(c)_i & = & \sum_{j\ge i}c_j &\\
         & = & max(\sum_{j\ge i}a_j,\sum_{j\ge i}b_j)& \mbox{by~Theorem~\ref{th_lbi}}\\
         & = & max(\pi(a)_i,\pi(b)_i) &\\
         & = & inf(\pi(a),\pi(b))_i &\mbox{in the Young lattice}
\end{array}
$$
which proves the claim.
\end{proof}

Notice that this order embedding is not a \emph{lattice} embedding,
since it does not preserve the supremum. For instance, if $a=(2,2)$
and $b=(1,1,1)$, then $\pi(a) = (4,2)$, $\pi(b) = (3,2,1)$, and $c =
sup(a,b) = (2,1)$ in \tlbi\ but $\pi(c) = (3,1)$ and
$sup((4,2),(3,2,1)) = (3,2)$ in the Young lattice. There can be no
lattice embedding from $\tlbi$ to the Young lattice since the fact
that this one is a \emph{distributive} lattice would imply that
$\tlbi$ would be distributive, which is not true. Finally, notice
that a study similar to the one presented in this paper can be found
in \cite{Lat01} on another kind of integer partitions, namely
$b$-ary partitions. The Young lattice is a particular case of the
lattices introduced in this paper, and the reader interested in the
relations between \tlbi\ and the Young lattice should refer to it.

\section{The infinite binary tree $T_B(\infty)$}

As shown in our procedure to construct \lbnp\ from \lbn, each element $a$ of
$\lbnp $ is obtained from an element $a'$ of $\lbn$ by addition of one
grain: $a = \inc{\mbox{$a'$}}{i}$ for some integer $i$. We will now
represent this relation by a tree where $a \in \lbnp$
is the son of $a' \in \lbn$ if and only if
$a = \inc{\mbox{$a'$}}{i}$ and we label with $i$ the edge
$a' \longrightarrow a$ in this tree. We denote this tree by \tb.
The root of this tree is the empty partition $()$.
We will show two ways to find the partitions of a given integer $n$
in \tb, which will make it possible to give an efficient and simple
algorithm to enumerate them. Moreover,
the recursive structure of this tree will allow us to obtain a recursive
formula for the cardinality of \lbn\ and some special classes of partitions.

From the construction of \lbnp\ from \lbn, it follows that
the nodes of this tree are the elements of $\bigsqcup_{n\ge 0}{\lbn}$,
and that each node $a$ has at least one son, \inc{a}{1},
and one more if $a$ begins with
a slippery plateau of length $l$: the element \inc{a}{\ell +1}.
Therefore, \tb\ is a binary tree. We will call
\emph{left son} the first of two sons, and
\emph{right son} the other (if it exists).
We call \emph{the level $n$ of the tree} the set of elements of depth $n$.
The first levels of \tb\ are shown in Figure~\ref{fig_arbre}.

\fig{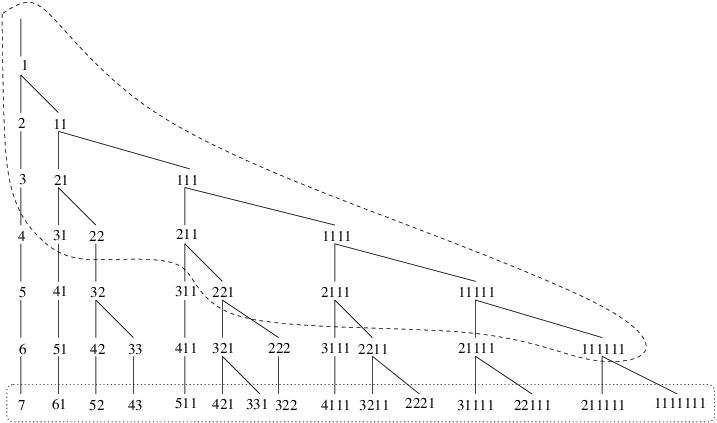}{The first levels of the tree $T_B(\infty)$ (to clarify
the picture, the labels are omitted). As shown on this
figure for $n=7$, we will discuss two ways to find the elements of
\lbn\ in \tb\ for any $n$.}

\noindent
Like in the case of \lbi, there are two ways to find
the elements of \lbn\ in \tb.
From the construction of \lbnp\ from \lbn\ given above, it is
straightforward that:
\begin{prop}
The level $n$ of \tb\ is exactly the set of the elements of \lbn.
\end{prop}

\noindent
Moreover, it is obvious from the construction of \tb\ that
the elements of the set
\mbox{$\overline{\lbnp} \setminus \overline{\lbn}$} are
sons of elements of $\overline{\lbn}$, therefore we deduce the following
proposition which can easily be proved by induction:
\begin{prop}
Let $\chi^{-1}$ be the inverse of the lattice isomorphism defined in
Lemma~\ref{th_chi}. Then, the set $\chi^{-1}(\overline{\lbn})$ is
a subtree of \tb\ having the same root.
\end{prop}

This proposition makes it possible to give a simple and efficient
algorithm to enumerate all the partitions of a given integer $n$,
using the binary tree structure it gives to the set of all these
partitions: Algorithm~\ref{algo_gener} acheives this in
linear time and space with respect to their number, which is optimal.

\begin{algorithm}
\begin{flushleft}
\SetVline
\In{An integer $n$}
\Out{The partitions of $n$}
\Begin{
 $\mbox{Resu} \leftarrow \emptyset$\;
 $\mbox{CurrentLevel} \leftarrow \lbrace () \rbrace$\;
 $\mbox{OldLevel} \leftarrow \emptyset$\;
 $l \leftarrow 0$\;
 \While{$\mbox{CurrentLevel} \not= \emptyset$}{
  \ForEach{$e$ in CurrentLevel}{
   Compute $p$ such that $p_i = e_{i-1}$ for all $i>1$ and $p_1 = n - l$\;
   Add $p$ to Resu\;
   }
  $\mbox{OldLevel} \leftarrow \mbox{CurrentLevel}$\;
  $\mbox{CurrentLevel} \leftarrow \emptyset$\;
  $l \leftarrow l + 1$\;
  \ForEach{$p$ in OldLevel}{
   Add $\inc{p}{1}$ to CurrentLevel\;
   \If{$p$ begins with a slippery plateau of length $l$}{
    Add $\inc{p}{\ell +1}$ to CurrentLevel\;
   }
  }
 \ForEach{$p$ in CurrentLevel}{
  \If{$n-l < p_1$}{
   Remove $p$ from CurrentLevel\;
   }
  }
 }
 Return(Resu)\;
}
\end{flushleft}
\caption{\label{algo_gener}Efficient computation of the partitions of an integer.}
\end{algorithm}

We will now give a recursive description of \tb.
We first define a certain kind of subtrees of \tb.
Afterwards, we show how the whole
structure of \tb\ can be described in terms of such subtrees.

\begin{definition}
We will call \emph{$X_k$ subtree} any subtree $T$ of\ \tb\ which is rooted
at an element \mbox{$a =(\lgr{i,\dots,i}{k},a_{k+1},\dots)$} with
$a_{k+1} \le i-1$ and which is either the whole subtree of \tb\ rooted at $a$ in the case $a$
has only one son, or $a$ and its left subtree otherwise.
Moreover, we define $X_0$ as a simple node.
\end{definition}

\noindent
The next proposition
shows that all the $X_k$ subtrees are isomorphic.

\begin{prop}
A $X_k$ subtree, with $k\ge 1$, is composed by a chain of $k+1$
nodes (the rightmost chain) whose edges are labeled by $1$, $2$,
$\dots$, $k$ and whose $i$-th node is the root of a $X_{i-1}$
subtree for all i between $1$ and $k+1$. (See
Figure~\ref{fig_def_Xk}.)
\label{prop_rec_xk}
\end{prop}
\begin{proof}
The claim is obvious for $k=1$. Indeed, in this case
the root $a$ has the form
$(i,a_2,\dots)$ with $a_2 \le i-1$, therefore its left son has the form
$(i+1,i-1,\dots)$, \ie\ it starts with a cliff, and has only one son.
This son also starts with a cliff; we can then deduce that $X_1$ is
simply a chain, which is the claim for $k=1$.

Suppose now that the claim is proved for any $i<k$ and consider the root $a$
of a $X_k$ subtree: \mbox{$a=(\lgr{i,\dots,i}{k},a_{k+1},\dots)$}
with $a_{k+1}\le i-1$. Its left son is \mbox{$\inc{a}{1} =
(i+1,i,\dots,i,a_{k+1},\dots)$} with $a_{k+1}\le i-1$, therefore it
is the root of a $X_1$ subtree. Moreover, \inc{a}{1}\ has one right
son: \mbox{$a^{\fb{1}\fb{2}} = (i+1,i+1,i,\dots,i,a_{k+1},\dots)$},
which by definition is the root of a $X_2$ subtree. After $k-1$ such
stages, we obtain \mbox{$a^{\fb{1}\fb{2}\dots\fb{k-1}}$}, which is
equal to \mbox{$(i+1,\dots,i+1,i,a_{k+1})$}. This node is the root
of a $X_{k-1}$ subtree and has a right son:
$$a^{\fb{1}\fb{2}\dots\fb{k-1}\fb{k}}, \ie\
(\lgr{i+1,\dots,i+1}{k},a_{k+1},\dots)$$
and we still have
$a_{k+1}\le i-1$. Therefore, this node is the root of a $X_k$
subtree, and from the definition of $T_B(\infty)$ we know that it
has no other son. This completes the proof.
\end{proof}

\fig{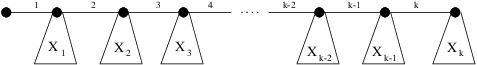}{Self-referencing structure of $X_k$ subtrees}

\noindent
This recursive structure and the above propositions allow us to
give a compact representation of the tree by a chain:

\begin{theorem}
The tree \tb\ can be represented by the infinite chain defined as
follows: the $i$-th node of this chain, $(\lgr{1,\dots,1}{i-1})$, is
linked to the following node in the chain by an edge labeled with
$i$ and is the root of a $X_{i-1}$ subtree. See
Figure~\ref{fig_arbre_chaine}. \label{th_arbrechaine}
\end{theorem}

\fig{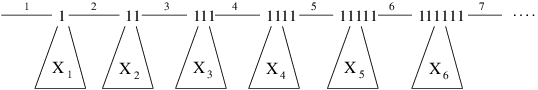}{Representation of $T_B(\infty)$ as a chain}

\noindent
Moreover, we can prove a stronger property of each subtree in this chain:

\begin{corol}
The $X_k$ subtree of \tb\ with root $(1,\dots,1)$
contains exactly the partitions of length $k$.
\label{th_part_born}
\end{corol}
\begin{proof}
Because of their recursive structure shown in
Proposition~\ref{prop_rec_xk},
$X_k$ subtrees contain no edge with label greater
than $k$. Therefore, if the root of a $X_k$ subtree is of length $k$ then
all its nodes have length $k$.
Moreover, no $X_l$ subtrees with $l \not= k$ and with a root of length $l$
can contain any node of length $k$.
This remark, together with Theorem~\ref{th_arbrechaine}, implies the result.
\end{proof}

\noindent
We can now state our last result:

\begin{corol}
Let $c(\ell,k)$ denote the number of paths in a $X_k$ tree originating from the
root and having length $\ell$. We have:
$$ c(\ell,k) = \left
\{\begin{array}{llll}
1 & \mbox{ if }  \ell = 0 \mbox{ or } k = 1\\
\sum_{i=1}^{inf(\ell,k)} c(\ell-i,i) & \mbox{ otherwise }
\end{array}
\right.$$
Moreover, $|L_B(n)| = c(n,n)$ and the number of partitions of $n$ with length
exactly $k$ is \mbox{$c(n-k,k)$}.
\end{corol}
\begin{proof}
The formula for $c(\ell,k)$ is derived directly from the structure of $X_k$ trees
(Proposition~\ref{prop_rec_xk} and Figure~\ref{fig_def_Xk}).
To obtain the formula for $|\lbn|$, we recall 
Proposition~\ref{prop_rec_xk} and Theorem~\ref{th_arbrechaine}. Applying
them, we observe that if we keep the nodes
of depth $n$ at most, then the subtrees obtained from \tb\ and $X_n$ turn
out to be isomorphic.
The last formula is directly
derived from Theorems~\ref{th_arbrechaine} and \ref{th_part_born}.
\end{proof}

\bibliographystyle{elsart-num-sort}
\bibliography{BIBLIO}

\end{document}